\documentclass{amsart}
\usepackage{amsmath}
\usepackage{amssymb}

    \newtheorem{theorem}{Theorem}
\numberwithin{theorem}{section} \theoremstyle{plain}

\newtheorem{question}[theorem]{Question}

\newtheorem{corollary}[theorem]{Corollary}

\newtheorem{proposition}[theorem]{Proposition}

\theoremstyle{definition}
\newtheorem{definition}[theorem]{Definition}
\newtheorem{notation}[theorem]{Notation}

\newtheorem{remark}[theorem]{Remark}
\numberwithin{equation}{section}

\newcommand{\fq}{\mathbb{F}_q}\newcommand{\cale}{\mathcal{E}}\newcommand{\kappap}{\kappa_P}
\newcommand{\qp}{q_P}\newcommand{\ep}{\mathcal{E}_P}\newcommand{\tep}{t(\ep)}
\newcommand{\thetaep}{\theta(\ep)}\newcommand{\bfe}{\mathbf{E}}
\DeclareMathOperator*{\eend}{End}

\begin{document}
\title[Distribution of the of Frobenius elements]
{On the distribution of the of Frobenius elements on elliptic curves over function fields}
\author{Am\'{\i}lcar Pacheco}
\date{\today}
\address{Universidade Federal do Rio de Janeiro (Universidade do Brasil)\\ Departamento de
Ma\-te\-m\'a\-ti\-ca Pura\\
Rua Guai\-aquil 83, Cachambi, 20785-050 Rio de Janeiro, RJ, Brasil}
\thanks{This work was partially supported by CNPq research grant 300896/91-3
and Pronex \#41.96.0830.00.\\ \indent Details can be found on \cite{pac2}} \email{amilcar@impa.br} \maketitle

\section{Introduction}
Let $C$ be a smooth projective irreducible curve defined over a finite field $\fq$ of $q$ elements and characteristic
$p>3$ with function field $K=\fq(C)$. Let $E/K$ be a non-constant elliptic curve and $\varphi:\cale\to C$ its minimal
regular model. For every $P\in C$, denote by $\deg(P)=[\kappap:\fq]$ its degree and $\qp=q^{\deg(P)}$. If
$\ep=\varphi^{-1}(P)$ is an elliptic curve over $\kappap$ denote by $\tep=\qp+1-\#\ep(\kappap)$ the \textit{trace of
Frobenius of $\cale$ at $P$}. It follows from a theorem of Hasse-Weil \cite[Chapter V, Theorem 2.4]{sil}
$\tep=\qp^{1/2}(e^{i\thetaep}+e^{-i\thetaep})=2\qp^{1/2}\cos(\thetaep)$ with $0\le\thetaep\le\pi$ and $\tep\le
2\qp^{1/2}$. Denote by $C_0$ the set of points $P\in C$ such that $\ep$ is an elliptic curve. In \cite{pac2} we
addressed the following question.

\begin{question}\label{ques1}
Let $B\ge 1$ and $t$ be integers with $|t|\le 2q^{B/2}$. Let $\pi(B,t)=\#\{P\in C_0\,|\,\deg(P)\le
B\text{ and }\tep=t\}$. How big is $\pi(B,t)$?
\end{question}

This question is an analogue for elliptic curves over function fields of the Lang-Trotter conjecture (cf. \cite{lt} and
\cite{mu1}). The goal of this note is to deduce from the answer of Question \ref{ques1} (cf. Corollary \ref{cora}) the
following more precise version of the Sato-Tate conjecture \cite{ta} (cf. Corollary \ref{corb}), in the sense that we
actually have an estimate for every $B\ge 1$, not an asymptotic result. The Sato-Tate conjecture was proved in the case
of elliptic curves over function fields by Yoshida \cite{yo} and Murty \cite{mu2}.

\begin{question}Let $B\ge 1$ be an integer and $\alpha,\beta\in[0,\pi]$ and $\pi(B,\alpha,\beta)=\#\{P\in
C_0\,|\,\deg(P)\le B\text{ and }\alpha\le\thetaep\le\beta\}$. How big is $\pi(B,\alpha,\beta)$?
\end{question}

For each $k\le B$ such that $|t|\le 2q^{k/2}$ we start by estimating $\pi(k,t)'=\#\{P\in C_0\,|\,\deg(P)=k\text{ and
}t(\mathcal{E}_P)=t\}$. Let $\mathcal{E}_P'=\mathcal{E}_P\times_{\kappa_P}\mathbb{F}_{q^k}$ and $\pi(k,t)''=\#\{P\in
C_0(\mathbb{F}_{q^k})\,|\,t(\mathcal{E}_P')=t\}$. The former set is contained in the latter one so
$\pi(k,t)'\le\pi(k,t)''$ and throughout all this paper we actually estimate $\pi(k,t)''$.

\section{Preliminaries}

Observe first that $E/K$ has to be an ordinary elliptic curve, otherwise $j(E)\in\mathbb{F}_{p^2}$ (cf. \cite[Chapter
V, Theorem 3.1]{sil}) , but this contradicts the fact that $E/K$ is non-constant.

Let $j_{\mathcal{E}}:C\to\mathbb{P}^1$ be the $j$-map induced from $\varphi_{\mathcal{E}}$. We say that $P\in C_0$ is
good ordinary, respectively good supersingular, if $\mathcal{E}_P$ is an ordinary, respectively supersingular, elliptic
curve. Since the number of supersingular $j$-invariants in $\overline{\mathbb{F}}_q$ is finite (cf. \cite[Chapter V,
Theorem 4.1]{sil}), then the number of good supersingular points $P\in C_0$ is also finite and bounded by an absolute
constant. This does not hold for elliptic curves over $\mathbb{Q}$ (cf. \cite{elk}). Recall that an elliptic curve
$\bfe/\mathbb{F}_q$ with $t(\bfe)=q+1-\#\bfe(\mathbb{F}_q)$ is supersingular if and only if $p\nmid t(\bfe)$ (cf.
\cite[Ex. 5.10]{sil}).

\section{Estimate of $\pi(k,t)''$}

\begin{notation}Let $I(t)$ be set of the isogeny classes of elliptic curves
${E}/\mathbb{F}_{q^k}$ defined over such that $\#{E}(\mathbb{F}_{q^k})=q^k+1-t$. Let $\mathfrak{A}_{k,t}$ be the set of
$\mathbb{F}_{q^k}$-isomorphism classes $[{E}]$ of ${E}\in I(t)$ and $N(t)=\#\mathfrak{A}_{k,t}$.
\end{notation}

\begin{definition}Let $\Delta<0$ be an integer such that $\Delta\equiv 0\text{ or }1\pmod 4$,
$B(\Delta)=\{\alpha x^2+\beta xy+\gamma y^2\,|\,\alpha,\beta,\gamma\in\mathbb{Z},\alpha>0\text{ and
}\beta^2-4\alpha\gamma=\Delta\}$ and $b(\Delta)=\{\alpha x^2+\beta xy+\gamma y^2\in
B(\Delta)\,|\,\gcd(\alpha,\beta,\gamma)=1\}$. The group $\text{SL}_2(\mathbb{Z})$ acts on $B(\Delta)$ via
$(\begin{smallmatrix}\alpha &\beta\\ \gamma &\delta\end{smallmatrix})f(x,y)=f(\alpha x+\beta y,\gamma x+\delta y)$
preserving $b(\Delta)$. The sets $b(\Delta)/\text{SL}_2(\mathbb{Z})$ and $B(\Delta)/\text{SL}_2(\mathbb{Z})$ are finite
with cardinality $h(\Delta)$ and $H(\Delta)$, respectively. The numbers $h(\Delta)$ and $H(\Delta)$ are called the
class number  and the Kronecker class number of $\Delta$,  respectively.
\end{definition}

\begin{proposition}\cite[Proposition 2.2]{sch}\label{prop1}
Let $\Delta<0$ be an integer such that $\Delta\equiv 0\text{ or }1\pmod 4$ then
\begin{equation}\label{eqn1}
H(\Delta)=\sum_fh(\Delta/f^2),
\end{equation}
where $f$ runs through all positive divisors of $\Delta$ such that $\Delta/f^2\in\mathbb{Z}$ and $\Delta/f^2\equiv
0\text{ or }1\pmod 4$.
\end{proposition}

\begin{remark}\label{rem1}
Let $\mathcal{O}$ be an imaginary quadratic order with discriminant $\Delta(\mathcal{O})$ and $h_{\mathcal{O}}$ its
class number. It follows from the correspondence between binary quadratic forms and complex quadratic orders that
$h_{\mathcal{O}}=h(\Delta(\mathcal{O}))$, where $\Delta(\mathcal{O})$ denotes the discriminant of $\mathcal{O}$
\cite[Chap. 2, Section 7, Theorem 4]{bs}.
\end{remark}

\begin{proposition}\cite[Theorem 4.5]{sch}
\label{prop2}Suppose $t\nmid p$ and let $\bfe\in I(t)$ with $\mathcal{O}=\eend_{\mathbb{F}_{q^k}}(\bfe)$.
$\#\{[\bfe']\in\mathfrak{A}_{k,t}\,| \,\mathcal{O}=\eend_{\mathbb{F}_q}(\bfe')\}=h_{\mathcal{O}}$.
\end{proposition}

\begin{notation}Denote $\mathcal{O}(t^2-4q^k)$ the imaginary quadratic order
 of discriminant $t^2-4q^k$.
 \end{notation}

\begin{proposition}\label{cor1}\cite[Theorem 4.6]{sch}
\begin{equation*}N(t)=\left\{\begin{aligned}H(t^2-4q^k),&\text{ if }t^2<4q^k\text{ and }p\nmid t\\
H(-4p),&\text{ if }t=0\text{ and }q^k\text{ is not a square}\\
\frac 1{12}\left(p+6-4\left(\frac{-3}p\right)-3\left(\frac{-4}p\right)\right),&\text{ if }t^2=4q^k\text{ and }q^k\text{
is a square}\\
1-\left(\frac{-3}p\right),&\text{ if }t^2=q^k\text{ and }q^k\text{ is a square}\\
1-\left(\frac{-4}p\right),&\text{ if }t=0\text{ and }q^k\text{ is a square}\\
0,&\text{ otherwise},\end{aligned}\right.\end{equation*} where $(\frac{\cdot}p)$ denotes the Legendre symbol at $p$.
\end{proposition}

\begin{proof}[Sketch of proof]We just prove the case where $p\nmid t$. In this case, it follows from
\cite[Theorem 4.3]{sch} that all imaginary quadratic orders $\mathcal{O}\supset\mathcal{O}(t^2-4q^k)$ occur as
$\mathbb{F}_{q^k}$-endomorphism ring of elliptic curves in $I(t)$. Hence, the result is a consequence of Propositions
\ref{prop1} and \ref{prop2} and Remark \ref{rem1}.
\end{proof}

\begin{theorem}\label{thm1}
$$\pi(k,t)''\le\deg(j_{\mathcal{E}})N(t).$$
\end{theorem}

\begin{proof}Let ${C}_0(\mathbb{F}_{q^k})$ be the set of
$\mathbb{F}_{q^k}$-rational points of $C_0$ and $\mathcal{C}_{k,t}=\{P\in
C_0(\mathbb{F}_{q^k})\,|\,t(\mathcal{E}_P')=t\}$, where
$\mathcal{E}_{P}'=\mathcal{E}_P\times_{\kappa_P}\mathbb{F}_{q^k}$. Define $\psi:\mathcal{C}_{k,t}\to\mathfrak{A}_{k,t}$
by $\psi(P)=[\mathcal{E}_{P}']$ and let $j(\mathcal{E}_P)$ be the $j$-invariant of $\mathcal{E}_P$.

We claim that $\psi^{-1}([\mathcal{E}_{P}'])\subset j_{\mathcal{E}}^{-1}(j(\mathcal{E}_{P}))$. In fact, if
$Q\in\psi^{-1}([\mathcal{E}_{P}'])$, then there exists an $\mathbb{F}_{q^k}$-isomorphism between $\mathcal{E}_{Q}'$ and
$\mathcal{E}_{P}'$, in particular $j(\mathcal{E}_{Q})=j(\mathcal{E}_{P})$. Hence,
$\#\psi^{-1}([\mathcal{E}_{P}'])\le\#j_{\mathcal{E}}^{-1}(j(\mathcal{E}_{P})) \le\deg(j_{\mathcal{E}})$ and
\begin{equation}\begin{split}\label{eqn0}\pi(k,t)''&=\sum_{[E]\in\psi(\mathcal{C}_{k,t})} \#\psi^{-1}([E])
\\ &\le\deg(j_{\mathcal{E}})\#\psi(\mathcal{C}_{k,t})\le\deg(j_{\mathcal{E}})
N(t).\end{split}\end{equation}
\end{proof}

\begin{corollary}\label{cora}
$$\pi(B,t)\le\left(\sum_{\substack{k\le B\\ |t|\le
2q^k}}N(t)\right)\deg(j_{\mathcal{E}}).$$
\end{corollary}

We now turn to the distribution of the angles $\thetaep$. Note that given an elliptic curve $\bfe/\mathbb{F}_{q^k}$, we
have $0\le\alpha\le\theta(\bfe)\le\beta\le\pi$ if and only if $t_{\beta}(k)\le t(\bfe)\le t_{\alpha}(k)$, where
$t_{\beta}(k)=2q^k\cos(\beta)$ and $t_{\alpha}(k)=2q^k\cos(\alpha)$. So, $\pi(B,\alpha,\beta)=\#\{P\in
C_0\,|\,\deg(P)=k\le B\text{ and }t_{\beta}(k)\le\tep\le t_{\alpha}(k)\}$.

\begin{corollary}\label{corb}
$$\pi(B,\alpha,\beta)\le\left(\sum_{\substack{k\le B\\ t_{\beta}(k)\le t\le t_{\alpha}(k)\\ t^2\le
4q^k}}N(t)\right)\deg(j_{\cale}).$$
\end{corollary}

\begin{remark}\label{remkey}We would like to compute examples in which we could test whether the
bound of Theorem \ref{thm1} is achieved. One good sort of example comes from modular curves. However, in this case
there is almost no control on $j_{\mathcal{E}}$ in contrast with the $j$-map $J$ naturally associated to the modular
problems. Moreover, if we observe the proof of Theorem \ref{thm1} closely, we notice that we can replace the regular
minimal model by any elliptic surface $\varphi_{\widetilde{\mathcal{E}}}:\widetilde{\mathcal{E}}\to C$ having ${E}/K$
as the generic fiber, defining the notions of good ordinary (good supersigular) points in terms of the fibers of
$\widetilde{\mathcal{E}}\to C$ being smooth ordinary (supersingular) elliptic curves. In this set-up it makes sense to
consider the trace of Frobenius of the fibers of good ordinary points. We can also consider elliptic curves
$\mathbb{E}\to C_1$ in the sense of \cite[Chapter 2]{km} defined over an affine subcurve $C_1\subset C$ with generic
fiber ${E}/K$ and compute the number (still denoted $\pi(k,t)''$) of $\mathbb{F}_{q^k}$-rational points $P\in C_1$
corresponding to good ordinary fibers $\mathbb{E}_P$ such that $t(\mathbb{E}_P')=t$. We will also denote by $\pi(B,t)$,
respectively $\pi(B,\alpha,\beta)$, the number of points $P\in C_1$ with $\deg(P)\le B$ and such that
$t(\mathbb{E}_P)=t$, respectively $\alpha\le\theta(\mathbb{E}_p)\le\beta$. The elliptic curve comes equipped with a
$j$-map $J:C\to\mathbb{P}^1$ and we look for conditions for $\pi(k,t)''$ to be equal to $\deg_s(J)N(t)$, where
$\deg_s(J)$ denotes the separable degree of $J$.
\end{remark}

We give explicit estimates of $\pi(B,t)$ and $\pi(B,\alpha,\beta)$ in the case of universal elliptic curves (cf.
Corollaries 5.3, 5.4, 5.6, 5.7, 5.9 and 5.10).

\section{Affine models}\label{affine}

Let $X$ be a smooth irreducible projective curve over $\mathbb{F}_q$ and $Y\subset X$ an affine subcurve. Suppose there
exists an elliptic curve $\mathbb{E}\to Y$ with generic fiber ${E}/K$ with $K=\fq(X)$ and a map $J:X\to \mathbb{P}^1$
whose restriction to $Y$ is given by $y\mapsto j(\mathbb{E}_y)$, where $\mathbb{E}_y$ denotes the fiber of
$\mathbb{E}\to Y$ at $y$. Given $y\in Y(\mathbb{F}_{q^k})$, let $\kappa_y$ be its residue field and
$\mathbb{E}'_y=\mathbb{E}_y\times_{\kappa_y}\mathbb{F}_{q^k}$. Let $\mathcal{Y}_{k,t}=\{y\in
Y(\mathbb{F}_{q^k})\,|\,t(\mathbb{E}_y')=t\}$ and $\pi(k,t)''=\#\mathcal{Y}_{k,t}$. Let
$\vartheta:\mathcal{Y}_{k,t}\to\mathfrak{A}_{k,t}$ be the map defined by $y\mapsto [\mathbb{E}_y']$.

\begin{proposition}\label{propkey}\cite[Proposition 3.1]{pac2}
Suppose the following three conditions are satisfied:
\begin{enumerate}
\item $\vartheta^{-1}([\mathbb{E}_{y}'])=J^{-1}(j(\mathbb{E}_{y}))$. \item $\vartheta$ is surjective. \item For every
$y\in Y$, the inertia degree $f(y\,|\,j(\mathbb{E}_y))$ equals $1$. The set $\mathcal{R}\subset Y$ of possible
ramification points of $J$ is contained in $J^{-1}(\{0,1728\})$. For each $y\in\mathcal{R}$ the ramification index
$e(y\,|\,0)$ (respectively $e(y\,|\,1728)$) of $P$ over $0$, respectively $1728$, equals $3$, respectively $2$.
\end{enumerate}
Then $\pi(k,t)''=\deg_s(J)H(t^2-4q^k)$, where $\deg_s(J)$ denotes the separable degree of $J$.
\end{proposition}

\section{Universal elliptic curves}

\subsection{Igusa curves}

Let $\bfe$ be an elliptic curve defined over a field $L$ of characteristic $p$. The absolute Frobenius $F_{\text{abs}}$
induces an isogeny $F_{\text{abs}}:\bfe\to\bfe^{(p)}$, where $\bfe^{(p)}$ denotes the elliptic curve obtained by
raising the coefficients of a Weierstrass equation of $\bfe$ to the $p$-th power. For each $n\ge 1$, let
$F^n_{\text{abs}}:\bfe\to\bfe^{(p^n)}$ be the $n$-th iterate of $F^n_{\text{abs}}$. Let $V^n$ be the dual isogeny of
the $n$-th iterate $F^n_{\text{abs}}$ of $F_{\text{abs}}$. An Igusa structure of level $p^n$ in $\bfe$ is a generator
of $\ker(V^n)$.

There exists a smooth affine curve $Y_n$ over $\mathbb{F}_p$ parametrizing isomorphism classes of pairs $(\bfe,P)$,
where $\bfe$ is an elliptic curve defined over an $\mathbb{F}_p$-scheme $S$ and $P\in\bfe^{(p^n)}(S)$ is an Igusa
structure of level $p^n$. In fact, $Y_n$ is a coarse moduli scheme for the moduli problem
$[\text{Ig}(p^n)]:\bfe/S/\mathbb{F}_p\mapsto P$. The compactification $X_n$ of $Y_n$ obtained by adding $\phi(p^n)/2$
points at infinity (called the cusps) is a smooth projective irreducible curve over $\mathbb{F}_p$ called the Igusa
curve of level $p^n$ \cite[Chapter 12]{km}, where $\phi$ denotes the Euler function.

An elliptic curve is $\bfe/S/\mathbb{F}_p$ is ordinary, if each of its geometric fibers is ordinary. An Igusa ordinary
(respectively Igusa supersingular) point $y\in Y_n$ is a point representing the isomorphism of a pair $(\bfe,P)$, where
$\bfe/S/\mathbb{F}_p$ is an ordinary elliptic curve, $S$ an $\mathbb{F}_p$-scheme (respectively $\bfe/L$ is a
supersingular elliptic curve, $L$ a field of characteristic $p$) and $P\in\bfe^{(p^n)}(S)$ (respectively
$P\in\bfe^{(p^n)}(L)$) is an Igusa structure of level $p^n$ in $\bfe$.

The group $(\mathbb{Z}/p^n\mathbb{Z})^*$ acts on $Y_n$ by $a\mapsto(\bfe,aP)$ and the group $\{\pm 1\}$ acts trivially.
These actions are extended to $X_n$ permuting the cusps simply transitively. Let $y\in Y_n$ represent the isomorphism
class of a pair $(\bfe,P)$. If $y$ is Igusa supersingular, then $y$ is fixed by $(\mathbb{Z}/p^n\mathbb{Z})^*$. If $y$
is Igusa ordinary and $j(\bfe)=1728$, respectively $j(\bfe)=0$, then $y$ has a stabilizer of order $2$, respectively
$3$, in $(\mathbb{Z}/p^n\mathbb{Z})^*/\{\pm 1\}$. On all other points of $Y_n$, $(\mathbb{Z}/p^n\mathbb{Z})^*/\{\pm
1\}$ acts freely. We identify the quotient of $X_n$ by $(\mathbb{Z}/p^n\mathbb{Z})^*/\{\pm 1\}$ to the projective line
$\mathbb{P}^1$ and the quotient map $J_n:X_n\to\mathbb{P}^1$ is Galois of degree $\phi(p^n)/2$. Its restriction to
$Y_n$ is given by $(\bfe,P)\mapsto j(\bfe)$.

The curve $Y_n^{\text{ord}}$ obtained from $Y_n$ by removing the Igusa supersingular points is a fine moduli space for
the restriction of $[\text{Ig}(p^n)]$ to ordinary elliptic curves. This means that there exists a universal elliptic
curve $\mathbb{E}_n\to Y^{\text{ord}}_n$ such that every ordinary elliptic curve $\bfe/S/\mathbb{F}_p$ with an Igusa
structure $P\in\bfe^{(p^n)}(S)$ of level $p^n$ is obtained from $\mathbb{E}_n\to Y^{\text{ord}}_n$ by a unique base
extension. In particular, if $K_n$ is the function field of $X_n$ over $\mathbb{F}_p$ and ${E}_n/K_n$ is the generic
fiber of $\mathbb{E}_n\to Y^{\text{ord}}_n$, then ${E}_n/K_n$ is the unique elliptic curve defined over $K_n$ with
$j$-invariant $j({E}_n)$ and a $K_n$-rational Igusa structure of level $p^n$.

If ${E}\in I(t)$ and $P\in{E}^{(p^n)}(\mathbb{F}_{p^k})$ is an Igusa structure of level $p^n$. Since ${E}$ and
${E}^{(p^n)}$ are isogenous, then $t\equiv p^k+1\pmod{p^n}$. So for the rest of this subsection, we assume $t\equiv
p^k+1\pmod{p^n}$.

\begin{proposition}\label{propa}\cite[Proposition 4.1]{pac2}
Conditions (1), (2) and (3) of Proposition \ref{propkey} are satisfied, a fortiori
$\pi(k,t)''=(\phi(p^n)/2)H(t^2-4p^k)$.
\end{proposition}

\begin{proof}In the notation of Section \ref{affine}, $Y=Y_n^{\text{ord}}$.
Let $y\in\mathcal{Y}_{k,t}$, denote $\mathbb{E}_{n,y}/\kappa_y$ the fiber of $\mathbb{E}_n\to Y_n^{\text{ord}}$ at $y$
and $\mathbb{E}_{n,y}'=\mathbb{E}_{n,y}\times_{\kappa_y}\mathbb{F}_{p^k}$.

(1) Let $x\in\vartheta^{-1}([\mathbb{E}_{n,y}'])$, then $\mathbb{E}_{n,x}'$ is $\mathbb{F}_{p^k}$-isomorphic to
$\mathbb{E}_{n,y}'$, in particular $j(\mathbb{E}_{n,x})=j(\mathbb{E}_{n,y})$, i.e., $x\in
J_n^{-1}(j(\mathbb{E}_{n,y}))$. Let $x\in J_n^{-1}(j(\mathbb{E}_{n,y}))$, then $x$ represents the isomorphism class of
the pair $(\mathbb{E}_{n,y},P_y)$, where $P_y\in\mathbb{E}_{n,y}^{(p^n)}(\kappa_y)$ is an Igusa structure of level
$p^n$. By the geometric description of $J_n$, there is no inertia at $Y_n^{\text{ord}}$, hence $\kappa_x=\kappa_y$.
Furthermore, $\mathbb{E}_{n,y}$ is an elliptic curve over $\kappa_y=\kappa_x$ with $j$-invariant equal to
$j(\mathbb{E}_{n,y})=j(\mathbb{E}_{n,x})$ having a $\kappa_y$-rational Igusa structure. It follows from the universal
property of $\mathbb{E}_n\to Y_n^{\text{ord}}$ that $\mathbb{E}_{n,y}=\mathbb{E}_{n,x}$, a fortiori
$[\mathbb{E}_{n,x}']=[\mathbb{E}_{n,y}']$ and $x\in\vartheta^{-1}([\mathbb{E}_{n,y}'])$.

(2) For every $[{E}]\in\mathfrak{A}_{k,t}$, $\#{E}(\mathbb{F}_{p^k})=\#{E}^{(p^n)}(\mathbb{F}_{p^k})\equiv
0\pmod{p^n}$. Thus, there exists an Igusa $P\in{E}^{(p^n)}(\mathbb{F}_{p^k})$ structure of level $p^n$. Let $y\in
Y_n^{\text{ord}}(\mathbb{F}_{p^k})$ represent the isomorphism class of the pair $({E},P)$. But $\mathbb{E}_{n,y}$ is
the unique elliptic curve over $\kappa_y$ with $j$-invariant $j(\mathbb{E}_{n,y})=j(E)$ having a $\kappa_y$-rational
Igusa structure of level $p^n$. Thus, $\mathbb{E}_{n,y}'={E}$. In particular, $[\mathbb{E}_{n,y}']=[{E}]$ and
$\vartheta$ is surjective.

Condition (3) follows from the geometric description of $J_n$. Consequently,
$$\pi(k,t)''=\frac{\phi(p^n)}2\sum_{\mathcal{O}(t^2-4p^k)\subset\,\mathcal{O}}
h_w(\Delta(\mathcal{O}))=\frac{\phi(p^n)}2H(t^2-4p^k).$$\end{proof}

\begin{remark}
Proposition \ref{propa} was implicitly used in \cite[Corollary 2.13]{pac1} to obtain an explicit expression for
$\#X_n(\mathbb{F}_{p^k})$.
\end{remark}

\begin{corollary}
$$\pi(B,t)\le\left(\sum_{\substack{k\le B\\ t^2<4p^k\\ t\equiv p^k+1\pmod{p^n}}}H(t^2-4p^k)\right)\frac{\phi(p^n)}2.$$
\end{corollary}

\begin{corollary}
$$\pi(B,\alpha,\beta)\le\left(\sum_{\substack{k\le B\\ t_{\beta}(k)\le t\le t_{\alpha}(k)\\ t^2<4p^k\\  t\equiv
p^k+1\pmod{p^n}}}H(t^2-4p^k)\right)\frac{\phi(p^n)}2.$$
\end{corollary}

\subsection{The modular curve $X(N)$}

Let $N>2$ be an integer not divisible by $p$. Let $\zeta\in\overline{\mathbb{F}}_p$ be a primitive $N$-th root of unity
and $\mathbb{F}_q=\mathbb{F}_p(\zeta)$. Let $Y(N)$ be the affine smooth curve defined over $\mathbb{F}_q$ parametrizing
isomorphism classes of triples $(\bfe,P,Q)$, where $\bfe$ is an elliptic curve defined over an $\mathbb{F}_q$-scheme
$S$ and $P,Q\in\bfe[N](S)$ is a Drinfeld basis for $\bfe[N](S)$ and $e_N(P,Q)=\zeta$ \cite[3.1]{km}, where $e_N$
denotes the $N$-th Weil pairing (cf. \cite[2.8]{km} and \cite[III, \S8]{sil}). In fact, it is a fine moduli space for
the modular problem $[\Gamma(N)]:{E}/S/\mathbb{F}_q\mapsto (P,Q)$ such that $e_N(P,Q)=\zeta$. The compactification
$X(N)$ of $Y(N)$ obtained by adding the cusps is a smooth projective irreducible curve defined over $\mathbb{F}_q$
\cite[Theorem 3.7.1]{km}.

The group $\text{SL}_2(\mathbb{Z}/N\mathbb{Z})$ acts on $Y(N)$ by $(\begin{smallmatrix}a &b\\ c&
d\end{smallmatrix})(\bfe,P,Q)\mapsto (\bfe,aP+bQ,cP+dQ)$ and the group $\{\pm 1\}$ acts trivially. If $y\in Y(N)$
represents the isomorphism class of a triple $(\bfe,P,Q)$, then $y$ has a stabilizer of order $3$, respectively $2$, if
$j(\bfe)=0$, respectively $j(\bfe)=1728$. On all other points of $Y(N)$, $\text{SL}_2(\mathbb{Z}/N\mathbb{Z})/\{\pm
1\}$ acts freely. The stabilizer at every cusp has order $N$ \cite[Theorem 6]{ig}. So we identify the quotient of
$X(N)$ by $\text{SL}_2(\mathbb{Z}/N\mathbb{Z})/\{\pm 1\}$ to the projective line $\mathbb{P}^1$. Let $\mathbb{E}(N)\to
Y(N)$ be the universal elliptic curve of $Y(N)$ and ${E}(N)/K(N)$ its generic fiber. The quotient map
$J(N):X(N)\to\mathbb{P}^1$ is Galois of degree $(1/2)N\phi(N)\psi(N)$ and its restriction to $Y(N)$ is given by
$({E},P,Q)\mapsto j(\bfe)$, where $\psi(N)=N\prod_{\ell\mid N}(1+(1/\ell))$ and $\ell$ runs over the divisors of $N$.

Denote $\mathcal{O}((t^2-4q^k)/N^2)$ the imaginary quadratic order of discriminant $(t^2-4q^k)/N^2$. Let $\bfe\in
I(t)$, we have $\bfe[N]\subset\bfe(\mathbb{F}_{q^k})$ if and only if $t\equiv q^k+1\pmod{N^2}$, $q^k\equiv 1\pmod N$
and $\mathcal{O}((t^2-4q^k)/N^2)\subset\eend_{\mathbb{F}_{q^k}}(\bfe)$ \cite[Proposition 3.7]{sch}.  Assume till the
end of this subsection that $t\equiv q^k+1\pmod{N^2}$ and $q^k\equiv 1\pmod N$. Let
$\mathfrak{A}_{k,t}'=\{[\bfe]\in\mathfrak{A}_{k,t}\,|\,
\mathcal{O}((t^2-4q^k)/N^2)\subset\eend_{\mathbb{F}_{p^k}}({E})\}$. By \cite[Theorem 4.9]{sch},
$\#\mathfrak{A}_{k,t}'=H((t^2-4q^k)/N^2)$. Note that $H((t^2-4q^k)/N^2)<H(t^2-4q^k)$.

\begin{proposition}\cite[Proposition 4.3]{pac2} Condition (1) and (3) of Proposition \ref{propkey} are satisfied.
However, $\pi(k,t)''=(1/2)N\phi(N)\psi(N)H((t^2-4q^k)/N^2)$.
\end{proposition}

\begin{corollary}
$$\pi(B,t)\le\left(\sum_{\substack{k\le B\\ t^2\le 4q^k\\ t\equiv q^k+1\pmod{N^2}\\ q^k\equiv 1\pmod N}}
H\left(\frac{t^2-4q^k}{N^2}\right)\right)\frac{N\phi(N)\psi(N)}2.$$
\end{corollary}

\begin{corollary}
$$\pi(B,\alpha,\beta)\le\left(\sum_{\substack{
k\le B\\ t_{\beta}(k)\le t\le t_{\alpha}(k)\\ t^2\le 4q^k\\ t\equiv q^k+1\pmod{N^2}\\ q^k\equiv 1\pmod N}}
H\left(\frac{t^2-4q^k}{N^2}\right)\right)\frac{N\phi(N)\psi(N)}2.$$
\end{corollary}

\subsection{The modular curve $X_1(N)$}

Let $N>4$ be an integer not divisible by $2$, $3$ and $p$. Let $Y_1(N)$ be the smooth affine curve defined over
$\mathbb{F}_p$ parametrizing isomorphism classes of pairs $(\bfe,P)$, where $\bfe$ is an elliptic curve $\bfe$ defined
over an $\mathbb{F}_p$-scheme $S$ and $P\in\bfe(S)$ is a point  of exact order $N$ \cite[Chapter 3]{km}. In fact, it is
a fine moduli space for the moduli problem $[\Gamma_1(N)]$ defined by $(\bfe/S/\mathbb{F}_p,P)\mapsto P$. The
compactification $X_1(N)$ of $Y_1(N)$ is a smooth irreducible projective curve defined over $\mathbb{F}_p$
\cite[Theorem 3.7.1]{km}.

Let $\mathbb{E}_1(N)\to Y_1(N)$ be the universal elliptic curve of $Y_1(N)$ and ${E}_1(N)/K_1(N)$ its generic fiber.
The $j$-map $J(N):X(N)\to\mathbb{P}^1$ factors through the Galois cover $X(N)\to X_1(N)$ of degree $N$, whose
restriction to $Y(N)$ maps to $Y_1(N)$ by $(\bfe,P,Q)\mapsto (\bfe,P)$. It induces the $j$-map
$J_1(N):X_1(N)\to\mathbb{P}^1$ whose restriction to $Y_1(N)$ is given by $(\bfe,P)\mapsto j(\bfe)$. Since $2,3\nmid N$,
given $y\in Y(N)$ and $y_1\in Y_1(N)$ such that $J(N)(y)=J_1(N)(y_1)$ equals $0$, respectively $1728$, then the
ramification index $e(y\,|\,y_1)$ equals $1$. A fortiori, $e(y_1\,|\,0)=3$, respectively $e(y_1\,|\,1728)=2$. Note also
that since there exists no inertia in $Y(N)\to\mathbb{A}^1$, then the same holds for $Y_1(N)\to\mathbb{A}^1$, thus
Condition (3) of Proposition \ref{propkey} is satisfied.

Observe that if $\bfe\in I(t)$ has a point $P\in\bfe(\mathbb{F}_q)$ of exact order $N$, then
$N\mid\#\bfe(\mathbb{F}_{p^k})$. The converse holds if $N$ is a prime number. We assume till the end of this subsection
that $N=\ell$ is a prime number different from $2$, $3$ and $p$, and $t\equiv p^k+1\pmod{\ell}$.

\begin{proposition}\cite[Proposition 4.4]{pac2} Conditions (1), (2) and (3) of Proposition \ref{propkey} are
satisfied, a fortiori $\pi(k,t)''=(1/2)(\ell^2-1)H(t^2-4p^k)$.
\end{proposition}

\begin{corollary}
$$\pi(B,t)\le\left(\sum_{\substack{k\le B\\ t^2\le 4p^k\\ t\equiv p^k+1\pmod{\ell}}}H(t^2-4p^k)\right)\frac{\ell^2-1}2.$$
\end{corollary}

\begin{corollary}
$$\pi(B,\alpha,\beta)\le\left(\sum_{\substack{k\le B\\ t_{\beta}(k)\le t\le t_{\alpha}(k)\\ t^2\le 4p^k\\ t\equiv
p^k+1\pmod{\ell}}}H(t^2-4p^k)\right)\frac{\ell^2-1}2.$$
\end{corollary}

\end{document}